\def\Bbb#1{{\bf #1}}

\def\fnote#1{\footnote}
\def\blacksquare{\hbox{\vrule width 4pt height 4pt depth 0pt}}

\def\cwleftpar#1#2{\leftskip #1 \rightskip #2 plus 1fill}
\def\cwrightpar#1#2{\leftskip #1 plus 1fill \rightskip #2}
\def\cwcenterpar#1#2{\leftskip #1 plus 1fill \rightskip #2 plus 1fill}
\def\cwfullpar#1#2{\leftskip#1\rightskip#2}

\def\cwoutdent#1#2{\llap{\hbox to #1{#2 \hss}}\ignorespaces}
\def\cwparbegin#1#2#3#4#5{
	\ifcase #1 \cwleftpar{#2}{#3}
	\or \cwrightpar{#2}{#3}
	\or \cwcenterpar{#2}{#3}
	\else \cwfullpar{#2}{#3}\fi
	\ifcase #4 \baselineskip = 1.5\baselineskip
	\or \baselineskip = 2\baselineskip
	\or \baselineskip = 3\baselineskip
	\else \baselineskip = 1\baselineskip\fi
	\ifdim #5 > 0in \else \noindent \fi
	\noindent\ignorespaces}
\documentclass{article}
\begin{document}
\advance \vsize by -2\baselineskip
\def\makeheadline{
{\noindent \folio \par
\vskip \baselineskip }}

\vspace*{2ex}

\noindent {\Large On some generalizations of the Jacobi identity. II}

\vspace*{2ex}

\noindent {\it Bozhidar Zakhariev Iliev}
\fnote{0}{\noindent $^{\hbox{}}$Permanent address:
Laboratory of Mathematical Modeling in Physics,
Institute for Nuclear Research and \mbox{Nuclear} Energy,
Bulgarian Academy of Sciences,
Boul.\ Tzarigradsko chauss\'ee~72, 1784 Sofia, Bulgaria\\
\indent E-mail address: bozho@inrne.bas.bg\\
\indent URL: http://theo.inrne.bas.bg/$^\sim$bozho/}

\vspace*{3ex}

{\bf \noindent Published in:
Bulletin de la Societe des Sciences  et  des  Lettres  de Lodz,
vol.42, No. 28, Serie: Recherches sur  les  Deformations,
vol.14, No. 138, pp. 93-102, 1992.}\\[1ex]
\hphantom{\bf Published: }
http://www.arXiv.org e-Print archive No.~math.DG/0306019\\[2ex]

2003 PACS numbers: 02.40Vh\\[2ex]

\noindent
{\small
The \LaTeXe\ source file of this paper was produced by converting a
ChiWriter 3.16 source file into
ChiWriter 4.0 file and then converting the latter file into a
\LaTeX\ 2.09 source file, which was manually edited for correcting numerous
errors and for improving the appearance of the text.  As a result of this
procedure, some errors in the text may exist.
}\\[2ex]

	\begin{abstract}
On the basis of the generalizations of the Jacobi identity found by the
author some identities satisfied by the curvature and torsion of a covariant
differentiation are derived. A kind of the generalized covariant
differentiation is proposed and a method of finding some of satisfied by them
identities is given in the correspondence to the curvatures concerned.
Possible applications for some physical problems are pointed out.
	\end{abstract}\vspace{3ex}

 {\bf 1. Introduction}
\nopagebreak

\medskip
In the previous paper [5] we announced a kind of many-point generalizations of the widely used Jacobi identity. The purpose of the present investigation is to give more complete version of [5] including applications in identities satisfied by some derivation operators and the corresponding curvatures.

 Let {\it A} be an Abelian group and
$[\cdot,\cdot]:A\times  A\to A$. We say that the
operation [ , ] satisfies the $p$-th Jacobi identity, $p\ge 2, (cf. [5])$ if
for every $A_{i}\in A, i\in I\neq   $ and $i_{1},\ldots  i_{p}\in I$, we have
\[
 {\bigl(}[A_{i_{1}},[A_{i_{2}},[\ldots  ,[A_{i_{p-1}},A_{i_{p}}]\ldots
]]]{\bigr)}_{[i_1,[i_2,[\ldots,[i_{p-1},i_p]\ldots]]}
\]
\[
=p[A_{i_{1}},[A_{i_{2}},[\ldots  ,[A_{i_{p-1}},A_{i_{p}}]\ldots  ]]],
\quad  p\ge 2,\qquad (1.1)
\]
where the (multiple) bracket  operation
$[\cdot ,[\cdot \ldots  ,[\cdot ,\cdot ]\ldots  ]]$ on the indices is defined as
follows:

Let {\it A}  be Abelian group, the set I be not empty, $p$ and $q$ to be
integers, $2\le p\le q, r_{1},\ldots  r_{p}\in \{1,2,\ldots  ,q\}, r_{a}\neq
r_{b}$for $a\neq b, a,b=1,\ldots  ,p$, and to any $i_{1},\ldots  i_{q}\in I$
to correspond some $A_{i_{1}}\in A$. Denoting the opposite element of $A\in
A$ by -A, we for $p=2$ define
\[
(A_{ i_1\ldots i_q })_{ [i_{r_1},i_{r_2}] }
:=A_{i_1\ldots i_q}
- A_{ i_{\tau _2(1)} \ldots i_{\tau_2(q) } }
\]
 where the permutation $\tau _{2}:\{1,\ldots  ,q\}  \to \{1,\ldots  ,q\}$ is
such that $\tau _{2}(a):=a$ for $a\in \{1,\ldots  ,q\}\backslash
\{r_{1},r_{2}\}, \tau _{2}(r_{1}):=r_{2}$, and $\tau _{2}(r_{2}):=r_{1}$.

Let  the permutation $\tau _{p}:\{1,\ldots  ,q\}  \to \{1,\ldots  ,q\}$ be
defined by $\tau _{p}(a):=a$ for $a\in \{1,\ldots  ,q\}\backslash
\{r_{1},\ldots  ,r_{p}\}, \tau _{p}(r_{1}):=r_{2}, \tau
_{p}(r_{2}):=r_{3},\ldots  , \tau _{p}(r_{p-1}):=r_{p}$and $\tau
_{p}(r_{p}):=r_{1}$. Then for $3\le p\le q$ we define
\[
A_{i_{1}}\ldots i_q)
_{[i_{r_1},[i_{r_2},[\ldots [ i_{r_{p-1}},i_{r_p}]\ldots ] }
\]
\[
:=(
  A_{i_{1}\ldots i_q}
- A_{i_{\tau_p}(1)\ldots i_{\tau_p(q)} }
)_{[i_{r_2},[\ldots [i_{r_{p-1}}, i_{r_p}]\ldots ] },
\]
 where in the first row the square brackets are $p-1$ and in the second row
are $p-2$.

If  (1.1) holds for every $p\le k$ and $[-A,B]=-[A,B], A,B\in A$, then we
call $[\cdot,\cdot]$ antisymmetric of order k. In the case when {\it A} is a
ring and $[A,B]=[A,B]_{-}:=AB-BA$ is the commutator of A and $B$ the identity
(1.1) is valid for every $p\ge 2$. If for [,] we have $[-A,B]=-[A,B]$, then
it is antisymmetric of order $k=2,3,4$ iff it satisfies the first $k$ of the
following equalities (Sect. 4)
\[
{\bigl(}[A_{i_{1}},A_{i_{2}}]{\bigr)}_{<i_1,i_2>}=0,\qquad (1.2)
\]
\[
{\bigl(}[A_{i_{1}},[A_{i_{2}},A_{i_{3}}]]{\bigr)}
_{<i_1,i_2,i_3>}=0
\quad  ({\it Jacobi\ identity}),\qquad (1.3)
\]
\[
\bigl( [A_{i_{1}},[A_{i_{2}},[A_{i_{3}},A_{i_{4}}]]]
+ [A_{i_{1}},[A_{i_{4}},A_{i_{3}},A_{i_{2}}]]]{\bigr)}
_{<i_1,i_2,i_3,i_4>}=0,
 \qquad (1.4)
\]
where $i_{1},i_{2},i_{3},i_{4}\in I$ and a cyclic summation over the indices
included in $<>$ is performed, i.e. $(\ldots  )_{<i_{1\ldots >}}$means the
sum of all elements of {\it A} obtained from $(\ldots  )$ by a cyclic
permutation of the indices $i_{1},\ldots  $ .

  For $A_{i_{1}\ldots i_p}\in A,\ p=2,3,4$, we have the identities
\[
\bigl( ( A_{i_{1}\ldots i_p}
+ (-1)^{p} A_{i_1 i_p i_{p-1} \ldots i_2} )
	_{ [i_1,[\ldots,[i_{p-1},i_p]\ldots] }
\bigr)_{<i_1\ldots,i_p>} \equiv  0  \qquad (1.5)
\]
 which reduce to $(1.2)-(1.4)$ respectively, whenever
{\it A} is a ring and $[A,B]=AB-BA$.

The generalized Jacobi identities, and in particular $(1.2)-(1.5)$,
definitely can find applications in the theory of physical systems described
by curvature- and/or torsion-depending Hamiltonians (or Lagrangians)
[7,8,9,13], especially those containing higher derivatives [9]. This is
because the curvature and torsion are simple combinations of commutators of
some derivations (see, e.g., (2.1) below or Sec. 4), and - on the other side-
the commutators are exactly binary operations for which (some of) the
generalized Jacobi identities hold. At the level of classical Hamiltonian
mechanics, with or without constraints, there is also a place for possible
application of the mentioned identities. In fact, as it is well known, the
Poisson bracket is a classical analogue of the quantum commutator is [3]. In
our terminology, the Poisson bracket is at least a third order antisymmetric
operation, since it is antisymmetric and satisfies the Jacobi identity [3].
Hence, in principle, the generalized Jacobi identities will give a
relationship between the curvature and torsion corresponding to some
commutators (or, e.g. Poisson brackets). This relationship may turn out to be
useful for the physical theories, for instance, for finding conservation laws
or first integrals of the equations of motion. Examples of this kind consist
in using the (second) Bianchi identity (in the form of the vanishing of the
covariant divergence of the Einstein's tensor) for the derivation of
conservation laws in general relativity [14], and in using the classical
Jacobi identity for generating first integrals of the equations of motion
from those already found ones [3].

\medskip
\medskip
 {\bf 2. Some Identities for the Covariant Differentiation}

\medskip
Consider a manifold $M$ endowed with an affine connection (covariant
differentiation) $\nabla  [2, 3]$. Let $A, B, C$ and $D$ be vector fields on
M. The curvature $R(A,B)$ and torsion $T(A,B)$ operators are defined by
\[
 R(A,B):=[\nabla _{A},\nabla _{B}]_{-}-\nabla _{[A,B]_{-}},\qquad (2.1a)
\]
\[
T(A,B):=\nabla _{A}B-\nabla _{B}A-[A,B]_{-},\qquad (2.1b)
\]
 where $[\nabla _{A},\nabla _{B}]_{-}:=\nabla _{A}\circ \nabla _{B}-\nabla
_{B}\circ \nabla _{A}$is the commutator of $\nabla _{A}$and $\nabla _{B}$and
$[A,B]_{-}:=A\circ B-B\circ A$ is the one of A and B.

Since the algebra of derivations along vector fields is also a ring and $[ ,
]_{-}$is bilinear, (1.1) is valid for $[ , ]_{-}$ (with $A_{i_{a}}=\nabla
_{B_{a}}, B_{a}$ being vector fields) for every $p\ge 2$ and, as a
consequence of this, the equalities $(1.2)-(1.3)$ are also true for $[ ,
]_{-}$. Let us see what these equalities mean  now.

From (1.2), we find
$0=([\nabla _{A},\nabla _{B}]_{-})_{<A,B>}=(\nabla
_{A}\circ \nabla _{B}-\nabla _{B}\circ \nabla _{A})_{<A,B>}=(R(A,B)+\nabla
_{[A,B]_{-}})_{<A,B>}=(R(A,B))_{<A,B>}$ or
\[
 R(A,B)+R(B,A)=0\qquad (2.2)
\]
which expresses the usual skew-symmetry of the curvature operator [2, 3, 4].

The Jacobi identity (1.3) yields the second Bianchi identity [2,3,4]. In
fact, we have
\[
 0=([\nabla _{A},[\nabla _{B},\nabla _{C}]_{-}]_{-})_{<A,B,C>}
\]
\[
=(\nabla _{A}\circ [\nabla _{B},\nabla _{C}]_{-}-[\nabla _{B},\nabla
_{C}]_{-}\circ \nabla _{A})_{<A,B,C>}
\]
\[
={\bigl(}(\nabla _{A}R)(B,C)+ \nabla _{A}\circ \nabla _{[B,C]_{-}}+
R(\nabla _{A}B,C)
\]
\[
+ R(B,\nabla _{A}C) - \nabla _{[B,C]_{-}}\circ  \nabla _{A}
\bigr)_{<A,B,C>}
\]
\[
= {\bigl(}(\nabla _{A}R)(B,C) +R(A,T(C,B))+\nabla
_{[A,[B,C]_{-}}{\bigr)}_{<A,B,C>}
\]
from where, using (2.2) and
\[
(\nabla _{[A,[B,C]_{-}})_{<A,B,C>}=\nabla _{([A,[B,C]_{-}}=\nabla _{0}=0,
\]
which is a consequence of (1.3), we get
\[
 {\bigl(}(\nabla _{A}R)(B,C)+R(T(A,B),C){\bigr)}_{<A,B,C>}=0.\qquad
(2.3)
\]

 For $p\ge 4$ the equality (1.1) produces new identities satisfied by
the curvature and torsion operators (2.1). Of course, these identities are
not independent of the known ones $[3, 4]$ but nevertheless they are new.
Beneath we shall derive only the first of them corresponding to $p=4$. In
that case (1.4) reduces to
\[
{\bigl(}[\nabla _{A},[\nabla _{B},[\nabla _{C},\nabla _{D}]_{-}]_{-}]_{-}+
(B\leftarrow\joinrel \to D){\bigr)}_{<A,B,C,D>}=0,
\]
 where $+(B\leftarrow\joinrel \to D)$ means that we have to add terms
obtained from the preceding ones by changing the symbols $B$ and D. Using the
definition of $[ , ]_{-}$and (2.1a), after some easy manipulation, from the
last equality we obtain:
\[
\bigl\{ (\nabla _{A}\nabla _{B}R)(C,D) + (\nabla _{A}R)
\bigl(
(B,\nabla _{D}C-[D,C]_{-})+(C,\nabla _{B}D-\nabla _{D}B)-(D,\nabla
_{B}C){\bigr)}
\]
\[
 + R{\bigl(}(A,-\nabla _{B}\nabla _{C}D+\nabla _{C}\nabla _{D}B+\nabla
_{D}[B,C]_{-}+[B,[C,D]_{-}]_{-})
\]
\[
 +(\nabla _{A}B,\nabla _{D}C+[C,D]_{-})+(\nabla _{A}C,\nabla _{B}D){\bigr)}
+ (B\leftarrow\joinrel \to D)\bigr\} _{<A,B,C,D>}=0 .
\]

Writing here explicitly the terms $+(B\leftarrow\joinrel \to D)$ and using
(2.1), we get after a simple calculations:
\[
 \bigl\{ {\bigl(}(R(A,B)+\nabla _{[A,B]_{-}})R{\bigr)}(C,D)+(\nabla
_{A}R)(B,[C,D]_{-})+(\nabla _{A}R)([B,C]_{-},D)
\]
\[
+ R(A,R(C,B)D) + R(A,R(C,D)B) - R(A,T(B,[C,D]_{-}))
\]
\[
  -R(A,T([B,C]_{-},D)) + R(T(A,B),[C,D]_{-})\bigr\} _{<A,B,C,D>}=0.
\qquad (2.4)
\]

Since $A, B, C$ and $D$ are arbitrary vector fields, the sum of the terms
containing commutators in this identity must be zero, i.e. (2.4) reduces to
the following two identities:
\[
 \bigl\{ {\bigl(}(R(A,B)R \bigr)
(C,D)+R(R(A,B)C,D)+R(C,R(A,B)D){\bigr\}}_{<A,B,C,D>}=0 \qquad (2.5)
\]
 and
\[
 {\bigl\{}(\nabla _{[A,B]_{-}}R)(C,D) + (\nabla _{A}R)(B,[C,D]_{-}) +
(\nabla _{A}R)([B,C]_{-},D)
\]
\[
R(A,T(B,[C,D]_{-})) - R(A,T([B,C]_{-},D))
\]
\[
 + R(T(A,B),[C,D]_{-}){\bigr\}}_{<A,B,C,D>}=0\qquad (2.6)
\]

It is easy to see that (2.6) is equivalent to
\[
 {\bigl\{}{\bigl(}(\nabla _{[C,D]_{-}}R)(A,B)+R(T([C,D],A),B)
\bigr) _{<[C,D],A,B>}{\bigr\}}_{<A,B,C,D>}=0,\qquad (2.6^\prime )
\]
 and hence it
is evident that (2.6) is  a  corollary  to  the  second Bianchi identity
(2.3). Let us note that  if  we  have  grouped  the terms into the identity
preceding (2.4) in a different way we  might get directly (2.5) without the
"additional" identity  (2.6).  Analogous identities, in which only the
first derivatives  of  $R$  appear, may be obtained from (1.1) for $p\ge 5$,
but, because of their more complicated algebraic structure, we will not
derive them here.

 Let us see now how the identities (1.5), which are more "general" then
$(1.2)-(1.4)$, can be used for obtaining identities for the torsion and
curvature.

For $p=2$ from (1.5) we get
\[
0=((\nabla _{A}B)_{[A,B]})_{<A,B>}
= (\nabla _{A}B-\nabla _{B}A)_{<A,B>}
\]
\[
=(T(A,B)+[B,A])_{<A,B>}=(T(A,B))_{<A,B>},
\]
 i.e.,
\[
 T(A,B)+T(B,A)=0\qquad (2.7)
\]
 which expresses the well known $[2, 3]$ skew-symmetry of the torsion.

For $p=3$ the identity (1,5) reproduces  the first Bianchi identity $[3, 4]$.
In fact, it yields
\[
0=((\nabla _{A}\nabla _{B}C)_{[A,[B,C]]})_{<A,B,C>}
\]
\[
=(\nabla _{A}\nabla _{B}C-\nabla _{B}\nabla _{C}A-\nabla _{A}\nabla
_{C}B+\nabla _{C}\nabla _{B}A)_{<A,B,C>}
\]
\[
 =(R(C,B)A+\nabla _{A}(T(B,C))+T([C,B]_{-},A)+[[C,B]_{-},A]_{-}){ }
_{<A,B,C>}
\]
 which, by means of $(2.2), (2.1a)$ and (1.3), can be rewritten as
\[
 {\bigl(}R(A,B)C-(\nabla _{A}T)(B,C)-T(T(A,B),C)
\bigr)_{<A,B,C>}=0\qquad (2.8)
\]
As one may expect, (1.5) for $p=4$ generates a
new identity for the curvature and torsion of the connection. Since the
derivation of this identity is simple but pretty long enough, as that of
(2.5), we shall skip the calculations and present here only some of the
intermediate results:
\[
 0=\{(\nabla _{A}\nabla _{B}\nabla
_{C}D)_{[A,[B,[C,D]]]}+(B\leftrightarrow D)\}_{<A,B,C,D>}=
\]
\[
 =\cdot \cdot \cdot =\{(\nabla _{A}\nabla _{B}T)(C,D)+(\nabla _{A}T)(\nabla
_{B}(C,D))+T(\nabla _{A}\nabla _{B}(C,D))
\]
\[
 +\nabla _{A}\nabla _{B}[C,D]_{-}-((\nabla _{A}R)(B,C))(D)-((\nabla
_{A}R)(C,D))(B)
\]
\[
 +(R(\nabla _{A}(B,C)))(D)-(R(\nabla _{A}(C,D)))(B)-R(C,D)(\nabla _{A}B)
\]
\[
 -\nabla _{A}(\nabla _{[B,C]_{-}}D+\nabla _{[C,D]_{-}}B)+\nabla
_{[C,D]_{-}}(\nabla _{B}A)+(B\leftarrow\joinrel \to D)\}_{<A,B,C,D>}=\cdot
\cdot \cdot
\]
\[
 =\{(\nabla _{A}T)(\nabla _{B}(C,D)+\nabla _{D}(C,B))\}_{<A,B,C,D>},
\]
 where $\nabla _{A}(B,C):=(\nabla _{A}B,C)+(B,\nabla _{A}C)$. After some easy
calculations, from the last result, finally, we find
\[
 \{(R(A,B))(T(C,D))-(R(A,B)T)(C,D)-T(R(A,B)C,D)
\]
\[
  -T(C,R(A,B)D)\}_{<A,B,C,D>}=0\qquad (2.9)
\]
 A feature of this new identity
is that if the connection is curvature or torsion free, then in the both
cases it reduces to the trivial one: $0=0$.

\medskip
\medskip
 {\bf 3. A generalization of the Covariant Differentiation}

\medskip
Let there be given a family $\xi :=\{\xi _{a}: a\in \Lambda \neq   \}$ of
real vector bundles $\xi _{a}=(E_{a},\pi _{a},M), \pi _{a}:E_{a}\to M$ over a
differentiable manifold $M$ and the dimensions of $\xi _{a}$and $\xi _{b}$be
equal for every $a,b\in \Lambda $, i.e. $\pi ^{-1}_{a}(x)$ and $\pi
^{-1}_{b}(x), x\in M$ be isomorphic vector spaces. The set of $C^{k}$sections
of $\xi _{a}, k\ge 0$, will be denoted by Sec$^{k}(\xi _{a})$; Sec$(\xi
_{a})$ means the set of all sections over $\xi _{a}$.

 Let us define maps
\[
  ^{a,b}I_{x\to y}:\pi ^{-1}_{a}(x)\to \pi ^{-1}_{b}(y),
\quad  x,y\in M,\qquad (3.1)
\]
 which will be called transports; we suppose that:
\[
  ^{b,c}I_{y\to z}\circ ^{a,b}I_{x\to y}=^{a,c}I_{x\to z},
\quad a,b,c\in \Lambda ,\quad  x,y,z\in M\qquad (3.2)
\]
and
\[
 ^{a,a}I_{x\to x}=id_{\pi ^{-1_{(x)}}_{a}}.\qquad (3.3)
\]
 One can easily prove

{\bf Proposition} ${\bf 3}{\bf .}{\bf 1}{\bf .} A$ map (3.1) satisfies (3.2)
and (3.3) if and only if there exist a vector space $Q$ of dimension
$\dim(\pi ^{-1}_{a}(x))$ and 1:1 maps $^{a}F_{x}:\pi ^{-1}_{a}(x)\to Q$ such
that
\[
  ^{a,b}I_{x\to y}=^{b}F^{-1}_{y}\circ ^{a}F_{x}.\qquad (3.4)
\]
 Further we shall need the maps $^{a,b}I_{x}$:Sec$(\xi _{a})\to $Sec$(\xi
_{a})$ defined by
\[
  (^{a,b}I_{x}(T))(y)
:=^{a,b}I_{x\to y}T(x),\quad T\in Sec(\xi _{a}),
\ x,y\in M.   \qquad (3.5)
\]
 Due to (3.2), we, evidently, have $^{b,c}I_{x}\circ
^{a,b}I_{y}=^{a,c}I_{y}$.

 On $\xi $ we can define the following generalization of the (linear)
covariant differentiation on vector bundles. Let us consider the maps
\[
  ^{a,b}\nabla _{V}:Sec^{1}(\xi _{a})\to Sec(\xi _{b}),\qquad (3.6)
\]
 $V$ being a vector field on $M$ which may be called a {\it formal
connection} on $\xi $. By definition the maps possess the properties
\[
  ^{a,b}\nabla _{V}(S+T)=^{a,b}\nabla _{V}S+^{a,b}\nabla _{V}T,\qquad (3.7)
\]
\[
^{a,b}\nabla _{V+W}=^{a,b}\nabla _{V}+^{a,b}\nabla _{W},\qquad (3.8)
\]
\[
 ^{a,b}\nabla _{f\cdot V}=f(x)\cdot ^{a,b}\nabla _{V},\qquad (3.9)
\]
\[
^{a,b}\nabla _{V}\circ (f\cdot )=^{a,b}I_{x}\circ (V(f)\cdot )+f(x)\cdot
^{a,b}\nabla _{V},\qquad (3.10)
\]
 where $S,T\in $Sec$^{1}(\xi _{a}), V$ and
$W$ are vector fields on $M, f:M\to {\Bbb R}, V(f)$ is the action of $V$ on
$f$ and $f\cdot $ means (left) multiplication with f.

It is  easily seen that the map
$\nabla :V\mapsto \nabla _{V}:\cup_{a\in\Lambda }Sec^{1}(\xi _{a})\to
\cup_{a\in\Lambda }Sec(\xi _{a})$ defined by $(\nabla
_{V}T)(x):=(^{a,a}\nabla _{V}T)(x)$ for every $T\in $Sec$^{1}(\xi _{a})$, has
the properties $\nabla _{V}(S+T)=\nabla _{V}S+\nabla _{V}T, \nabla
_{V+W}=\nabla _{V}+\nabla _{W}, \nabla _{f\cdot V}=f\cdot \nabla _{V}$and
$\nabla _{V}\circ (f\cdot )=V(f)+f\cdot \nabla _{V}$. Hence $\nabla $ defines
a covariant differentiation in every vector bundle $\xi _{a}[3, 4]$.

The conditions $(3.7)-(3.10)$ imply some restrictions on the used
transports. Namely, from (3.8) and (3.10) one immediately derives that
$^{a,b}I_{x\to y}((\sigma +\tau )T_{x})=^{a,b}I_{x\to y}(\sigma
T_{x})+^{a,b}I_{x\to y}(\tau T_{x})$ for every $\sigma ,\tau \in {\Bbb R}$
and $T_{x}\in \pi ^{-1}_{a}(x)$. In particular this means that the transports
must be ${\Bbb Z}$-linear. The condition is naturally satisfied if the
transports are ${\Bbb R}$-linear, i.e., if
\[
 ^{ _{a,b}}_{ }I_{x\to y}(\sigma S_{x}+\tau T_{x})=\sigma ^{a,b}I_{x\to
y}S_{x}+\tau ^{a,b}I_{x\to y}T_{x}, \sigma ,\tau \in {\Bbb R}\qquad (3.11)
\]
for every $S_{x},T_{x}\in \pi ^{-1}_{a}(x)$, which is assumed hereafter
everywhere.

The  transports $^{a,b}I_{x\to y}$and the derivatives $^{a,b}\nabla _{V}$will
be called {\it consistent} if
\[
 ^{a,b}\nabla _{V}=^{c,b}I_{y}\circ ^{a,c}\nabla _{V}.\qquad (3.12)
\]
 From this definition for $c=a$ immediately follows

 {\bf  Proposition 3.2.} If $^{a,b}I_{x\to y}$and $^{a,b}\nabla _{V}$are
consistent, then
\[
  ^{a,b}\nabla _{V}=^{a,b}I_{x}\circ \nabla _{V}\qquad (3.13)
\]
and, on the  opposite, if $\nabla _{V}:\cup_{\in\Lambda }Sec^{1}(\xi
_{a})\to \cup_{a\in\Lambda }Sec(\xi _{a})$ preserves the type of the
sections, i.e., $\nabla _{V}$:Sec$^{1}(\xi _{a})\to $Sec$(\xi _{a})$ and has
the written above properties, then the maps (3.13) satisfy $(3.7)-3.10)$, are
consistent with $^{a,b}I_{x\to y}$and the map Sec$^{1}(\xi _{a})\ni T\mapsto
(^{a,a}\nabla _{V}T)$ coincides with $\nabla _{V}$.

 If the maps $^{a,b}I_{x\to y}$are $C^{1}($with respect to $x$ or $y)$, we
can define
 \[
  ^{a,b}\nabla ^{I}_{V}
:=\lim_{\epsilon\to_0}
[\frac{1}{\varepsilon}(^{a,b}I_{x_{\epsilon }}-^{a,b}I_{x})]
={\bigl(} \frac{\partial}{\partial\varepsilon}
(^{a,b}I_{x_{\epsilon }}){\bigr)}\Big|_{\epsilon =0},\qquad (3.14)
\]
 where  $x^{\alpha }_{\epsilon }=x^{\alpha }+\epsilon V^{\alpha },
V=V^{\alpha }\partial /\partial x^{\alpha }$in some local basis.

 {\bf  Proposition 3.3.} The maps (3.14) satisfy $(3.7)-(3.10)$ and
\[
 ^{b,c}\nabla ^{I}_{V}\circ ^{a,b}I_{y}\equiv 0,\qquad (3.15)
\]
\[
^{b,c}\nabla ^{I}_{V}\circ ^{a,b}\nabla ^{I}_{W}\equiv 0,\qquad (3.16)
\]
\[
 ^{b,c}I_{y}\circ ^{a,b}\nabla ^{I}_{V}=^{a,c}\nabla ^{I}_{V},\qquad (3.17)
\]
the first of which does not depend on the linearity of $^{a,b}I_{x\to y}$.

{\it Proof}{\bf .} The first statement of the proposition follows directly
from $(3.14); (3.15)$ and the consistency condition (3.17) are consequences
from (3.14) and $^{b,c}I_{x}\circ ^{a,b}I_{y}=^{a,c}I_{y}; (3.16)$ is a
corollary of (3.17) and (3.15).\blacksquare

At the end of this section we shall write some expressions in local bases.

Let $\{e^{a}_{i}, i=1,\ldots  ,\dim(\xi _{a})\}$ and $\{\partial /\partial
x^{\alpha }\}$ be bases in Sec$(\xi _{a})$ and in the bundle tangent to $M$
respectively. We define the  transports
\[
 ^{a,b}I_{x\to y}e^{a}_{i}(x)=:^{a,b}H^{j}_{.i}(y,x)e^{b}_{j}(y),\qquad
(3.18)
\]
\[
 {\bigl(}(^{a,b}\nabla _{\partial /\partial x^{\alpha }})
(e^{a}_{i}){\bigr)}(y)=:^{a,b}\Gamma ^{j}_{.i\alpha }(y,x)e^{b}_{j}(y),
\qquad (3.19)
\] where the summation is understood to be
performed on the repeated on different levels indices within the range of
their values. One can easily verify that the consistency condition (3.12) is
equivalent to
\[
  ^{a,b}\Gamma ^{j}_{.i\alpha }(z,x)=^{b,c}H^{j}_{.k}(z,y)^{a,c}\Gamma
^{k}_{.i\alpha }(y,x)\qquad (3.20)
\]
 and if $T=T^{i}e^{a}_{i}$ and $V=V^{i}\partial /\partial x^{\alpha }$, then
\[
 ^{a,b}I_{x\to y}T(x)=^{a,b}H^{j}_{.i}(y,x)T^{i}(x)e^{b}_{j}(y),\qquad (3.21)
\]
\[
 (^{a,b}\nabla _{V}T)(y)=V^{\alpha }(x)[^{a,b}H^{j}_{.i}(y,x)(\partial
T^{i}/\partial x^{\alpha })(x)
\]
\[
+ ^{a,b}\Gamma ^{j}_{.i\alpha }(y,x)T^{i}(x)]e^{b}_{j}(y),\qquad (3.22)
\]
\[
(^{a,b}\nabla ^{I}_{V}T)(y)=V^{\alpha }(x)[^{a,b}H^{j}_{.i}(y,x)(\partial
T^{i}/\partial x^{\alpha })(x)+^{a,b}H^{j}_{.i\alpha
}(y,x)T^{i}(x)]e^{b}_{j}(y),\qquad (3.23)
\]
\[
((^{b,c}\nabla _{V}\circ ^{a,b}I_{y})(T))(z)=V^{\alpha }(x)
\bigl(-^{b,c}H^{j}_{.i\alpha }(z,x)
\]
\[
+  ^{b,c}\Gamma ^{j}_{.i\alpha }(z,x)
\bigr) ^{a,b}H^{i}_{.k}(x,y)T^{k}(y)e^{b}_{j}(z),\qquad (3.24)
\]
 where $(cf. (3.4))$
\[
^{a,b}H^{i}_{.j\alpha }(y,x)
:= \frac{\partial\,^{a,b}H_{.j}^{i}(y,x) }{\partial x^\alpha}
\]
\[
 =-^{a,b}H^{j}_{.k}(y,x)
	\frac{\partial\,^{b,a}H_{.l}^{k}(x,y) }{\partial x^\alpha}
{} ^{a,b}H^{l}_{.i}(y,x)\qquad (3.26)
\]
 From the equality (3.24) follows the evident implication
\[
  ^{b,c}\nabla _{V}\circ ^{a,b}I_{y}=0 \Leftrightarrow { }
^{b,c}\nabla _{V}=^{b,c}\nabla ^{I}_{V}.\qquad (3.27)
\]
 The purpose of what  follows is to find the analogues of curvature and
torsion operators for the derivatives $^{b,c}\nabla _{V}$and some satisfied
by them identities.

\medskip
\medskip
 {\bf 4. Curvatures for the Generalized Covariant Differentiation}

\medskip
Classically $[2, 4]$ the curvature and the torsion of a linear connection
arise from the consideration of combinations quadratic with respect to the
connection or some vector fields $(cf. (2.1))$, and these combinations are
skew-symmetric $(cf. (2.2))$. In our case a number of different such
expressions can be formed which may be divided into two groups. The first of
them contains the expressions obtained from
\[
 ^{b,c}\nabla _{W}\circ ^{a,b}\nabla _{V}\qquad (4.1)
\]
by antisymmetrizing it with respect to: $V$ and $W; x$ and $y; (x,V)$  and
$(y,W); b$ and $c; (b,V)$  and  $(c,W); (b,x)$  and  $(c,y); (b,x,V)$  and
$(c,y,W)$. The first three antisymmetrizations are formed in  a  usual way,
e.g. that one with respect to $V$ and $W$ gives
\[
  ^{b,c}\nabla _{W}\circ ^{a,b}\nabla _{V}-{ } ^{b,c}\nabla _{V}\circ
^{a,b}\nabla _{W}.
\]
 The remaining antisymmetrizations have two variants, e.g. that one with
respect to $(b,x)$ and $(c,y)$ gives
\[
  ^{b,c}\nabla _{W}\circ ^{a,b}\nabla _{V}-{ } ^{b,c}I_{z}\circ ^{c,b}\nabla
_{W}\circ ^{a,c}\nabla _{V},
\]
\[
 ^{b,c}\nabla _{W}\circ ^{a,b}\nabla _{V}-{ } ^{b,c}I\circ ^{c,b}\nabla
_{W}\circ ^{a,c}\nabla _{V},
\]
 where $^{b,c}I$:Sec$(\xi _{b})\to $Sec$(\xi _{c})$ and if $T\in $Sec$(\xi
_{b})$, then $(^{b,c}I(T))(x):=^{b,c}I_{x\to x}(T(x))$. Here the maps
$^{b,c}I_{z}$and $^{b,c}I$ appear because otherwise the differences are not
defined (the different terms in them map one and the same section into
sections from, generally, different bundles).

A  simple calculation, with the usage of (3.22) and (3.4), shows that if
$T\in $Sec$^{1}(\xi _{a})$ and $x,y,t\in M$, then
\[
 (^{b,c}\nabla _{W}\circ ^{a,b}\nabla _{V}(T))(t)=W^{\alpha }(y)V^{\beta
}(x)
\Big\{\Big[
^{b,c}H^{j}_{.i\alpha }(t,y)^{a,b}\Gamma ^{i}_{.k\beta }(y,x)
\]
\[
 ^{b,c}H^{j}_{.i}(t,y)
\Bigl( \frac{\partial {}^{a,b}\Gamma ^{i}_{.k\beta}(z,x)}
	    {\partial z^\alpha} \Big|_{z=y}
\Big] T^{k}(x)+[-^{b,c}H^{j}_{.i\alpha }(t,y)
\]
\[
 +^{b,c}\Gamma ^{j}_{.i\alpha }(t,y)]\delta ^{\gamma }_{\beta }((^{a,b}\nabla
_{\partial /\partial x^{\gamma }})(T))^{i}(y)
\Big\} e^{c}_{j}(y).\qquad (4.2)
\]
From here we see that any (first or higher) generalized covariant derivative
of some $C^{1}$section contains only first partial derivatives of that
section.

 The skewsymmetric expressions of the second class are obtained by
antisymmetrizing the expression
 \[
  ^{b,c}\nabla _{W}\circ ^{a,b}\nabla _{V},\qquad (4.3)
\]
 where $^{a,b}\nabla
:V\mapsto ^{a,b}\nabla _{V}$:Sec$^{1}(\xi _{a})\to $Sec$(\xi _{b})$ is such
that $(^{a,b}\nabla _{V}T)(x):=(^{a,b}\nabla _{V}T)(x)$ for $T\in
$Sec$^{1}(\xi _{a})$, with respect to $V$ and $W, b$ and $c$ and $(b,V)$ and
$(c,W)$.

  Applying (4.3) to $T\in $Sec$^{2}(\xi _{a})$ and using twice (3.22), we get
\[
  (^{b,c}\nabla _{W}\circ ^{a,b}\nabla _{V}(T))(x)=^{b,c}I_{x\to x}
\bigl(
((^{a,b}\nabla _{W(V^{\alpha _{)\partial /\partial x}}})T)(x){\bigr)}
\]
\[
+W^{\beta }(x)V^{\alpha }(x)\{[^{b,c}\Gamma ^{j}_{.i\beta
}(x,x)^{a,b}\Gamma ^{i}_{.k\alpha }(x,x) +{ } ^{b,c}H^{j}_{.i}(x,x)
\]
\[
\frac{d {}^{a,b}\Gamma_{.k\alpha}^{i}(x,x) }{d x^\beta}
-{ } ^{a,b,c}K^{\hbox{j.}\gamma }_{.l\beta \alpha }\Gamma
^{l}_{.k\gamma }(x,x)]T^{k}(x)+^{a,b,c}K^{\hbox{j.}\gamma }_{.k\beta \alpha }
\]
\[
 ((^{a,a}\nabla _{\partial /\partial x^{\gamma }})T)^{k}(x) +{ }
^{a,c}H^{j}_{.k}(x,x)
\frac{\partial^2 T^k(x) }
     {\partial x^\beta\partial x^\alpha}
\}e^{c}_{j}(x),\qquad (4.4)
\]
 where
\[
 ^{ _{a,b,c}}_{ }K^{\hbox{j.}\gamma }_{.k\beta \alpha }
:=^{b,c}H^{j}_{.i}(x,x)(
\frac{d {}^{a,b}H_{.k}^{l}(x,x) }
     {d x^\beta}
\delta ^{\gamma }_{\alpha }+{ } ^{a,b}\Gamma
^{i}_{.k\alpha }(x,x)\delta ^{\gamma }_{\beta })
\]
\[
  +{ } ^{b,c}\Gamma ^{j}_{.i\beta }(x,x)^{a,b}H^{i}_{.k}(x,x)\delta ^{\gamma
}_{\alpha },\qquad (4.5)
\]
 in which $\delta ^{\gamma }_{\alpha }=1$
for  $\alpha =\gamma $ and $\delta ^{\gamma }_{\alpha }=0$ for $\alpha \neq
\gamma $. The quantities (4.5) are symmetric with respect to $\alpha $ and
$\beta $ in the following important cases:
\[
^{a,a,a}K^{\hbox{j.}\gamma }_{.k\alpha \beta }=^{a,a}\Gamma
^{j}_{.k\alpha }(x,x)\delta ^{\gamma }_{\beta }+{ } ^{a,a}\Gamma
^{j}_{.k\beta }(x,x)\delta ^{\gamma }_{\alpha },\qquad (4.6a)
\]
\[
 + ^{a,b,c}K^{\hbox{j.}\gamma }_{.k\alpha \beta }|_{\nabla =\nabla
^{I}}=^{a,c}H^{j}_{.k\alpha }(x,x)\delta ^{\gamma }_{\beta }+{ }
^{a,c}H^{j}_{.k\beta }(x,x)\delta ^{\gamma }_{\alpha },\qquad (4.6b)
\]
 where
$\nabla =\nabla ^{I}$ means that $^{a,b}\nabla _{V}$ has to be replaced by
$^{a,b}\nabla ^{I}_{V}$the action of which on $T\in $Sec$^{1}(\xi _{a})$ is
$((^{a,b}\nabla ^{I}_{V})(T))(x):=(^{a,b}\nabla ^{I}_{V}T)(x)$.

 If we make the above-pointed antisymmetrizations in (4.2) and (4.4), which
are 11 and 5 in number respectively, we shall obtain expressions which are
linear with respect to $W, V, T$ and the first generalized covariant
derivatives of $T ($higher derivatives do not occur). The coefficients before
the last two quantities are, by definition, the components of the
corresponding curvatures. As the consideration of all of these 16 cases is
similar, we are going to consider, for brevity, only one of them.

The antisymmetrization of (4.3) with respect to $W$ and $V$, due to (4.4),
gives
\[
  ^{b,c}\nabla _{W}\circ ^{a,b}\nabla _{V}-{ } ^{b,c}\nabla _{V}\circ
^{a,b}\nabla _{W}-{ } ^{b,c}I\circ ^{a,b}\nabla _{[W,V]_{-}}
\]
\[
 =^{a,b,c}R(W,V) +{ } ^{a,b,c}D(W,V),\qquad (4.7)
\]
 where
\[
  (^{a,b,c}R(W,V)T)(x)=^{a,b,c}R^{j}_{.i\beta \alpha }T^{i}(x)W^{\beta
}(x)V^{\alpha }(x)e^{c}_{j}(x),\qquad (4.8a)
\]
\[
 (^{a,b,c}D(W,V)T)(x)=^{a,b,c}D^{\hbox{j.}\gamma }_{.i\beta \alpha }(((\nabla
_{\partial /\partial x^{\gamma }})T)(x))^{i}W^{\beta }(x)V^{\alpha
}(x)e^{c}_{j}(x), (4.8b)
\]
define the first and the second curvature operators  (for $^{a,b}\nabla )$
whose components are:
\[
 ^{a,b,c}R^{j}_{.i\beta \alpha }={\bigl(}^{b,c}H^{j}_{.k}(x,x)
\frac{d}{d x^\beta}
(^{a,b}\Gamma ^{k}_{.i\alpha }(x,x)) +{ } ^{b,c}\Gamma ^{j}_{.k\beta
}(x,x)
\]
\[
 ^{a,b}\Gamma ^{k}_{.i\alpha }(x,x){\bigr)}_{[\alpha ,\beta ]}-{ }
^{a,b,c}D^{\hbox{j.}\gamma }_{.k\beta \alpha }\Gamma ^{k}_{.i\gamma
}(x,x)\qquad (4.9a)
\]
\[
^{a,b,c}D^{\hbox{j.}\gamma }_{.i\beta \alpha
}:=(^{a,b,c}K^{\hbox{j.}\gamma }_{.i\beta \alpha })_{[\alpha ,\beta ]}.\qquad
(4.9b)
\]
 It is important to note that by virtue of (4.6), we have
\[
  ^{a,a,a}D^{\hbox{j.}\gamma }_{.i\beta \alpha }:=^{a,b,c}D^{\hbox{j.}\gamma
}_{.i\beta \alpha }|_{\nabla =\nabla ^{I}}=0.\qquad (4.10)
\]

 The first of
these cases, $b=c=a$, corresponds to the usual covariant differentiation in
the vector bundles $\xi _{a}$. In fact, in this case
$^{a,a,a}R(W,V)T=R(W,V)T=:^{a}R(W,V)T$, where $R(W,V)$ is given by (2.1a)
(see the above definition of $\nabla :V\mapsto \nabla _{V})$, so $^{a}R$ is
the usual curvature operator for $\nabla $ in $\xi _{a}$, the components of
which are obtained from (4.9a) for $c=b=a$, i.e. they are
\[
 ^{a}_{x}R^{j}_{.i\beta \alpha }
=\Bigl[
\frac{d}{d x^\beta}
 (^{a,a}\Gamma
^{j}_{.i\alpha }(x,x))+^{a,a}\Gamma ^{j}_{.k\beta }(x,x)^{a,a}\Gamma
^{k}_{.i\alpha }(x,x)
\Bigr]_{[\alpha,\beta ]}
\]
 (see (4.10) and (3.3)). Let us note that $^{a}_{x}R^{j}_{.i\beta \alpha }=0$
iff $\nabla =\nabla ^{I}$, i.e. iff $^{a,a}\Gamma ^{k}_{.i\alpha
}(x,x)=^{a,a}H^{k}_{.i\alpha }(x,x)$.

Now we are going to describe some of the identities satisfied by the above
curvatures.

Since the left-hand side of (4.7) can not be written in terms of commutators
of $^{a,b}\nabla _{V}$and $^{b,c}\nabla _{W}($they simply are not well
defined) some of the identities satisfied by the above curvatures may be
found on the base of (1.5) or from (1.1) for a special definition of the
operation [ , ]. Below, as examples, we consider the consequences of (1.5)
for $p=2,3$.

 In our case (1.5) for $p=2$ gives
\[
 0={\bigl(}(^{b,c}\nabla _{W}\circ ^{a,b}\nabla _{V})_{[W,V]}
\bigr) _{<W,V>}
=(^{b,c}\nabla _{W}\circ ^{a,b}\nabla _{V}-{ } ^{b,c}\nabla
_{V}\circ ^{a,b}\nabla _{W})_{<W,V>}
\]
\[
 =(^{b,c}\nabla _{W}\circ ^{a,b}\nabla _{V}-{ } ^{b,c}\nabla _{V}\circ
^{a,b}\nabla _{W}-{ } ^{b,c}I\circ ^{a,b}\nabla _{[W,V]_{-}})_{<W,V>}
\]
\[
=(^{a,b,c}R(W,V) +{ } ^{a,b,c}D(W,V))_{<W,V>}.
\]
Hence, because of (4.8), it follows that
\[
 ^{a,b,c}R(W,V)+^{a,b,c}R(V,W)=^{a,b,c}D(W,V)+^{a,b,c}D(V,W)=0\qquad (4.11)
\]
i.e. the skewsymmetry of the curvatures with respect to their vector
arguments.

 In the case under consideration (1.5) for $p=3$ gives:
 \[
 0={\bigl(}(^{c,d}\nabla _{A}\circ ^{b,c}\nabla _{B}\circ ^{a,b}\nabla
_{C})_{[A,[B,C]]}{\bigr)}_{<A,B,C>}
\]
\[
 ={\bigl(}(^{c,d}\nabla _{A}\circ ^{b,c}\nabla _{B}\circ ^{a,b}\nabla
_{C}-^{c,d}\nabla _{B}\circ ^{b,c}\nabla _{C}\circ ^{a,b}\nabla
_{A})_{[B,C]}{\bigr)}_{<A,B,C>}.
\]

 Therefore, taking into account (4.3), we get
\[
  [^{c,d}\nabla _{A}\circ (^{a,b,c}R(B,C)+^{a,b,c}D(B,C)+^{b,c}I\circ
^{a,b}\nabla _{[B,C]_{-}})
\]
\[
 -(^{b,c,d}R(B,C)+^{b,c,d}D(B,C)+^{c,d}I\circ ^{b,c}\nabla _{[B,C]_{-}})\circ
^{a,b}\nabla _{A}]_{<A,B,C>}=0.  \qquad (4.12)
\]
Because of (4.8) and (4.4), this
equality is equivalent to two identities which include the curvatures and
their first derivatives. These identities are obtained by equating to zero
the coefficients before $T$ and $(^{a,a}\nabla _{\partial /\partial x^{\gamma
}})T$ when (4.12) is applied to an arbitrary $T\in $Sec$^{2}(\xi _{a})$. The
identities pointed out, as well as their derivation, are simple but too long
to be presented here.

\medskip
\medskip
 {\bf 5. Conclusions}

\medskip
In this work we gave simple examples of application of the generalizations of the classical Jacobi identity announced in [5].

In the case of the covariant differentiation the first and the second of mentioned identities reproduce the known identities satisfied by the curvature and torsion tensors, but beginning with the third one, which is explicitly derived, the proposed here method generates new identities for these tensors.

In Sect. 3 we have proposed a kind of generalization of the covariant differentiation in a family of equidimensional vector bundles over a given differentiable manifold (which formally is
supposed to be endowed with a kind of a "transport along that manifold"). Here we wish to point out that a typical example of such families are the tensor bundles of a given rank. For instance, in the notions of Sect. 3, we may set $\xi =\{\xi _{1},\xi _{2},\xi _{3},\xi _{4}\}$ with $\xi _{1}=T(M)\otimes T(M), \xi _{2}=T(M)\otimes T^{*}(M), \xi _{3}=T^{*}(M)\otimes T(M)$ and $\xi _{4}=T^{*}(M)\otimes T^{*}(M)$, where $T(M)$ is the bundle tangent to $M$ and $T^{*}(M)$ is its dual. In section 8 we have described how the curvatures of this differentiation operations can be introduced and how on the basis of the generalized Jacobi identities some identities satisfied by them can be obtained. Because, in the general case, there are 16 such identities, we consider in details only two of them. An analogous treatment is still valid for the remaining identities.

\medskip
\medskip
 {\bf References}

\medskip
\medskip
[1] Chevalley C., The theory of Lie Groups, vol. I, Princeton Univ. Press, Princeton, $1946, ch$.IV, \S11.\par
[2] Faith C., Algebra: Rings, Modules and Categories, vol. I, Springer Verlag, $1973, ch. 1,3$.\par
[3] Goldstein $H$, Classical Mechanics, Addison-Wesley Press, 1951.\par
[4] Helgason S., Differential Geometry, Lie Groups and Symmetric Spaces, Academic Press, New York-San Francisco-London, 1978, p.69.\par
[5] Iliev B. Z. (1992), "On some generalizations of the Jacobi identity I.", Bull. Soc. Sci. Lettres L\'odz {\bf 42}, no.21, S\'erie: Recherches sur les d\'eformations {\bf 14}, no.$131, 5-11$, to appear.\par
[6] Kobayashi S., Nomizu K., Foundations of Differential geometry, vol I., Interscience Publishers, New York- London, 1963. \par
[7] Lanczos C., Electricity and general relativity, Rev. Mod. Phys., vol.{\bf 29}, No.$3, 1957, 337-350$.\par
[8] Lawrynowicz J. and  Succi F. (1993), "Almost supercomplex structures and special Fueter equations related to the many-particle problem", this volume, pp.\par
[9] Nesterenko V. V., "Singular Lagrangians with higher
derivatives", J. Phys. {\bf A}: Math. Gen., vol.{\bf 22}, No.$10, 1989, 1673-1687$.\par
[10] Sattinger D.H., O.L. Weaver, Lie groups and Algebras with Applications to Physics, Geometry and Mechanics, Springer Verlag, 1986.\par
[11] Schafer R.D., Introduction to Nonassociative Algebras, Academic Press, New York-London, 1963.\par
[12] Sternberg S., Lectures on Differential Geometry, Chelsea Publ. Co., New York, N. Y., 1983, p.360.\par
[13] Szczyrba W., A Hamiltonian structure of the interacting gravitational and matter fields, J. Math. Phys., vol.{\bf 22}, No.$9, 1981, 1926-1944$.\par
[14] Weinberg $S$, Gravitation and cosmology, John Wiley \& Sons, Inc., New
York-London-Sydney-Toronto, 1972.

\vspace{4ex}
 Institute for Nuclear Researchand Nuclear Energy

 Bulgarian Academy of Sciences

 Blvd. Tzarigradsko Chauss\'ee 72, 1784 Sofia, Bulgaria

\end{document}